\documentclass[11pt]{article}
\usepackage[a4paper,margin=1in]{geometry}
\usepackage{amsmath,amssymb,booktabs}
\usepackage[title]{appendix}
\usepackage[T1]{fontenc}
\usepackage{lmodern}
\usepackage{microtype}
\usepackage{hyperref}
\usepackage{float}
\usepackage{booktabs}

\hypersetup{colorlinks=true,linkcolor=blue,citecolor=blue,urlcolor=blue}

\newcommand{\dd}{\,d}

\title{Asymptotic expansions of the eigenvalues and norms for perturbations about the Falkner--Skan solution}
\author{Carlos Lozano* and Jorge Ponsin  \\
Theoretical and Computational Aerodynamics  \\
National Institute of Aerospace Technology (INTA), Spain  \\
*Corresponding author: lozanorc@inta.es }
\date{}

\begin{document}
\maketitle

\begin{abstract}

We derive large-mode asymptotic expansions for the eigenvalues and eigenfunction norms of the Chen-Libby perturbation problem about Falkner-Skan boundary-layer profiles. The analysis extends Brown's matched-asymptotic construction for the Blasius problem to Falkner-Skan solutions with positive wall shear---both favorable gradients and the upper branch of the mild adverse range. Although the outer and middle layers retain the same structure as in the Blasius case, the wall layer changes qualitatively when the pressure-gradient parameter $\beta$ is nonzero: the Falkner-Skan wall expansion singularly perturbs the inner Bessel problem at relative order $\Lambda^{-1/3}$, where $\Lambda$ is the large eigenvalue parameter. This produces a new wall-induced contribution $s^{-1/3}$ (where $s=n-1$ and $n$ is the large eigenvalue index) to the large-mode eigenvalue expansion, ahead of Brown's $s^{-1/2}$ correction. The new term vanishes in the Blasius limit, where Brown's ordering is recovered. The matched eigenfunction also yields asymptotic estimates for the Chen-Libby norms. The eigenvalue and norm formulae are compared with direct numerical shooting calculations for $\beta=1/2$, showing the expected asymptotic convergence and confirming the role of the new wall correction.

\end{abstract}

\section{Introduction}
Chen and Libby  \cite{ChenLibby1968} showed that small departures from a Falkner-Skan boundary layer are governed by a discrete eigenvalue problem analogous to the one that arises in the description of perturbations of the Blasius solution \cite{LibbyFox1963}. For adverse-pressure gradients, there are two solutions for the Falkner-Skan equation: one with positive friction coefficient (the upper branch) and another one with negative friction coefficient (the lower branch). For the upper branch, the resulting eigenvalue spectrum is positive and gives spatial decay of perturbations; on the lower branch, the spectrum is real but no longer sign-definite---both positive and negative eigenvalues occur~\cite{ChenLibby1968}---and the linearized problem is spatially unstable.  For favorable pressure gradients, the corresponding eigenvalues are positive and describe downstream relaxation to the similar Falkner-Skan profile.  


The relevance of the Chen-Libby eigenvalue problem has been highlighted by recent
studies on boundary layer sensitivity and non-modal stability. For instance, \cite{VaqueroRenard2019} showed that the
leading Chen-Libby eigenvalue governs the slow exponential decay of two-dimensional upstream
distortions, determining the long development lengths required for a perturbed laminar boundary
layer to recover its self-similar state. In the context of optimal streaks and reduced-order streak
dynamics in Falkner-Skan wedge flows, \cite{SanchezHigueraVega2011,HigueraVega2012} found that these eigenvalues form the exact limit of
the full three-dimensional spectrum as the spanwise wavenumber approaches zero. Because these
eigenvalues enter both the downstream relaxation of planar disturbances and the asymptotic
boundary of three-dimensional streak stability, obtaining precise asymptotic expansions provides
a direct way to understand their dependence on the pressure gradient parameter.

A central point in this perturbative boundary-layer analysis is the determination
of the eigenvalues. Since these arise by demanding exponential decay at infinity
rather than simple vanishing, the procedure is numerically cumbersome and has
generated a considerable literature \cite{Libby1970}. Soon after Chen-Libby's work,
Kotorynski \cite{Kotorynski1969}  put the Chen-Libby problem into Liouville normal form and
obtained a simple analytic approximation to the eigenvalues. His approximation
keeps the leading wall singularity and the leading harmonic growth of the
far-field potential, so that the resulting model is explicitly solvable in terms
of Laguerre polynomials and gives an equally spaced spectrum.

Brown  \cite{Brown1968}  developed a different, matched-asymptotic analysis for the Blasius
($\beta=0$) perturbation equation. Brown's analysis is based on the reduced
second-order equation first discussed by Stewartson \cite{Stewartson1957}. The eigenvalue
condition is obtained by constructing the solution in three overlapping regions,
matching them together, and imposing decay at infinity and regularity at the
wall. Kemp \cite{Kemp1970}, in the Blasius context, also showed how useful norm estimates
can be obtained by integrating an appropriately normalized approximate
eigenfunction. This observation motivates the norm calculation in Section~\ref{sec5},
although the eigenfunction used here is the matched Brown eigenfunction rather
than Kotorynski's Laguerre approximation.

A point that appears to have been overlooked in this classical problem is that
the Falkner-Skan generalization is not obtained from Brown's Blasius
calculation by a smooth replacement of constants. In the Blasius case, the
leading wall-layer problem is a Bessel equation whose large-argument phase
supplies the wall matching condition at the order required for Brown's first
correction. For a nonzero Falkner-Skan pressure-gradient parameter, however, the $\beta$-dependent terms of
the near-wall expansion combine to perturb the inner Bessel problem singularly at relative order
$\Lambda^{-1/3}$ in the wall scaling $\eta=O(\Lambda^{-1/3})$. As a consequence, the large-mode eigenvalue expansion acquires
a wall-induced $s^{-1/3}$ correction, which is absent in the Blasius limit and
precedes the usual Brown-type $s^{-1/2}$ term.

The purpose of the present note is to extend Brown's matched-asymptotic
calculation to the Chen-Libby problem for fixed favorable pressure-gradient
parameter $\beta>0$, and to use the resulting composite eigenfunction to
estimate the corresponding Chen-Libby norms. The calculation closely follows Brown's approach. 
After the reduction of the third-order equation
to a second-order equation for the transformed eigenfunction $z$, the solution is constructed in the same
three overlapping regions used by Brown: an outer turning-point region, a middle
WKB region, and a wall region. The large parameter is
\begin{equation}
        \Lambda=\lambda-\frac{1}{2}-2\beta,
        \label{eq:Lambda}
\end{equation}
and $s$ denotes the (large) mode index. A Liouville normal form is useful for
identifying the origin of the wall singularity but is not used below as a
separate eigenvalue approximation; the eigenvalue condition is obtained by
Brown's matching procedure. The same construction applies to the upper branch of
the adverse range $-0.1988<\beta<0$, where $a>0$ as well; this extension is
summarized in Section~\ref{sec6}.

The structure of the paper is as follows. In Section~\ref{sec2} we recall the
Chen-Libby eigenvalue problem and its reduction to Brown's second-order
equation. Section~\ref{sec3} develops the large-\(\Lambda\) asymptotics in the outer,
middle and wall regions. Section~\ref{sec4} carries out the matching and derives the
large-mode eigenvalue formula, including the constant that replaces Brown's
constant in the Falkner-Skan problem. Section~\ref{sec5} uses the same matched
eigenfunction to obtain an approximate formula for the
Chen-Libby norms, validated against direct numerical computation, and also discusses the Blasius limit.
Section \ref{sec6} briefly discusses the extension to adverse pressure gradients.

\section{Statement of the problem}
\label{sec2}

Perturbations about a Falkner-Skan boundary layer are described by a linear partial differential equation whose solution by means of separation of variables leads to the following eigenvalue problem \cite{ChenLibby1968}:
\begin{equation}
        N'''+fN''+(\lambda-2\beta)f'N'
        +(1-\lambda)f''N=0,
        \label{eq:CL}
\end{equation}
with
\begin{equation}
        N(0)=N'(0)=0,\qquad N'(\infty)=0
        \label{eq:CLbc}
\end{equation}
where an exponential decay of $N'$ at infinity is assumed in order to obtain a discrete eigenspectrum. The coefficients of \eqref{eq:CL} depend on the function $f(\eta;\beta)$ defined by
\begin{equation}
        f'''+ff''+\beta\{1-(f')^2\}=0,
        \qquad
        f(0)=f'(0)=0,
        \qquad
        f'(\infty)=1 
\label{eq:FS}
\end{equation}
which is the Falkner-Skan equation describing laminar viscous flow along a flat plate with outer velocity $u_\infty=Kx^\frac{\beta}{2-\beta}$ \cite{White}. The parameter $\beta$ can take both positive and negative values, corresponding to favorable and adverse pressure gradients, respectively. Solutions exist for $\beta\geq\beta_0 = -0.1988$. Two solution branches coexist for $\beta_0 < \beta < 0$, which coalesce at $\beta = \beta_0$. The upper branch describes forward flow (with $f'(\eta)>0$ for $\eta>0$ and positive wall shear stress $a=f''(0)>0$), while the lower branch describes unstable reverse flow with negative wall shear stress. For $\beta\geq 0$ only the upper branch survives. The case $\beta = 0$ corresponds to the Blasius boundary layer.  The
construction below is restricted to the case $a>0$ and $f'(\eta)>0$ for $\eta>0$. We develop it for
fixed favorable or zero gradient $\beta\geq0$ and extend it to the upper branch of the adverse range in Section~\ref{sec6}.

Equation \eqref{eq:CL} can be put in reduced form through the transformation \cite{Brown1968}

\begin{equation}
        N=f'U,
        \qquad
        U'=z\exp\left[-\frac12\int_0^\eta f(s)\,ds\right].
        \label{ztransf}
\end{equation}
The equation satisfied by \(z\) is
\begin{equation}
        f'z''+3f''z'+f'A_\beta(\eta;\Lambda)z=0,
        \label{eq:zEq}
\end{equation}
where
\begin{equation}
        \Lambda=\lambda-\frac12-2\beta,
\end{equation}
where
\begin{equation}
        A_\beta
        =
        \Lambda f'
        +\frac52\frac{f'''}{f'}
        -\frac14 f^2
        +\frac{\beta}{2}\left(f'-\frac{1}{f'}\right).
        \label{eq:A_beta}
\end{equation}
For \(\beta=0\), this reduces to Brown's equation, with \(\Lambda=\lambda-1/2\). The boundary conditions for $z$ can be obtained from \eqref{eq:CLbc} and \eqref{ztransf} as 
\begin{equation}
        az'(0)-\beta z(0)=0,
        \label{eq:zWallBC}
\end{equation}
at the wall and $z\to 0$ exponentially as $\eta\to\infty$. 

\section{Layer structure of the large \texorpdfstring{$\Lambda$}{Lambda} problem }
\label{sec3}

We now let $\Lambda$ be large. Three overlapping regions can be distinguished based on the behavior of the coefficient function $A_\beta$ in \eqref{eq:zEq}.
The outer region (large $\eta$) layer is centered around $\eta_0\sim\Lambda^{1/2}$. The wall ($\eta\to 0$) region has
$\eta=O(\Lambda^{-1/3})$ and leads to leading order to Brown's Bessel problem.
The middle region connects these two limits:
\begin{equation}
        \Lambda^{-1/3}\ll \eta,\qquad
        |\eta-\eta_0|\gg \Lambda^{-1/6}.
\end{equation}
In this region the coefficient of $z$ is positive and slowly varying on the
local oscillation scale. The eigenvalue condition is obtained by matching the
decaying outer solution to the middle WKB solution and then matching the latter
to the regular wall solution.
  
The calculation therefore proceeds as in Brown.  The appropriate form of equation \eqref{eq:zEq} in each layer is determined and the approximate solutions are obtained and matched at the corresponding overlap regions, yielding the desired condition on the eigenvalues. We now discuss each layer separately. 

\subsection{The outer approximation}
In this region, $\eta$ is large and the Falkner-Skan function may be replaced by its asymptotic form \cite{coppel}
\begin{equation}
        f(\eta)=\eta-\kappa_\beta+o(1),
        \qquad
        f'(\eta)=1+o(1),
        \label{largeeta}
\end{equation}
where
\begin{equation}
        \kappa_\beta=\int_0^\infty\{1-f'(\eta)\}\dd\eta .
\end{equation}
Thus
\begin{equation}
        A_\beta(\eta;\Lambda)
        =\Lambda-\frac14(\eta-\kappa_\beta)^2+o(1),
        \label{Aouter}
\end{equation}
which changes sign at the outer turning point
\begin{equation}
        \eta_0=\kappa_\beta+2\sqrt{\Lambda}+o(1).
        \label{outturnpnt}
\end{equation}
In view of \eqref{Aouter} we introduce the local variable 
\begin{equation}
        \xi=\frac{\eta-\kappa_\beta}{2\sqrt{\Lambda}}.
\end{equation}
Using \eqref{Aouter}, (\ref{eq:zEq}) reduces, up to algebraically small errors in $\Lambda$, to
\begin{equation}
        \frac{\dd^2z}{\dd\xi^2}+4\Lambda^2(1-\xi^2)z=0 .
        \label{eq:outereq}
\end{equation}
The coefficient of $z$ in \eqref{eq:outereq} has a zero at $\xi=1$.  For $\xi>1$ the coefficient is negative, and positive for $\xi<1$.  The standard Liouville-Green connection formulae \cite{JeffreysJeffreys,WhittakerWatson}, used by Brown in the 
Blasius problem, give the following solution for the eigenfunction $z$ in this regime as 
\begin{equation}
        z_{\rm out}^{\xi>1}(\xi)\sim \frac{1}{2} B_{\rm out}(\xi^2-1)^{-1/4}
        \exp\left\{-2\Lambda\int_1^\xi (y^2-1)^{1/2}\dd y\right\},
        \qquad \xi>1
\label{eq:outerdecay}
\end{equation}
which decays exponentially as $\eta\to\infty$, and
\begin{equation}
        z_{\rm out}^{\xi<1}(\xi)\sim B_{\rm out}(1-\xi^2)^{-1/4}\cos \Theta,
        \qquad 0<\xi<1.
\label{eq:outerosc}
\end{equation}
where
\begin{equation}
        \Theta=2\Lambda\int_\xi^1(1-y^2)^{1/2}\dd y-\frac{\pi}{4}.
        \label{zetaout}
\end{equation}

\subsection{The middle layer}

In the middle layer, $\eta=O(1)$ is the appropriate variable.  Dividing \eqref{eq:zEq} by $f'$ gives
\begin{equation}
        z''+3\frac{f''}{f'}z'
        +\left\{\Lambda f'+\frac{V_\beta}{f'}\right\}z=0,
        \label{eq:middleEq}
\end{equation}
where
\begin{equation}
        V_\beta=\frac{5}{2}f'''-\frac{1}{4}f^2f'
        +\frac{\beta}{2}\{(f')^2-1\}.
        \label{Vbeta}
\end{equation}

Let $k(\eta)=\sqrt{f'(\eta)}$ and
\begin{equation}
        t=\int_0^\eta\sqrt{f'(y)}\dd y.
        \label{eq:tauPhaseNorm}
\end{equation}
We write, following Brown,
\begin{equation}
        z=\Re\{g(\eta)e^{iY(\eta)}\},
        \qquad
        Y(\eta)=-x(\eta)+\gamma ,
        \label{eq:middleAnsatz}
\end{equation}
where
\begin{equation}
        x(\eta)=\sqrt{\Lambda} t.
        \label{eq:xdef}
\end{equation}
Substituting into \eqref{eq:middleEq} gives
\begin{equation}
        g''+
        \left(3\frac{f''}{f'}-2i\sqrt{\Lambda} k\right)g'
        +
        \left\{
        \frac{V_\beta}{f'}
        -i\sqrt{\Lambda}\left(k'+3\frac{f''}{f'}k\right)
        \right\}g=0 .
        \label{eq:geq}
\end{equation}
Let us expand $g$ in the form
\begin{equation}
        g=g_0-\frac{i}{\sqrt{\Lambda}}g_1+O(\Lambda^{-1})
\end{equation}
so that
\begin{equation}
        z_{\rm mid}
        =
        g_0\cos Y
        +
        \Lambda^{-1/2}g_1\sin Y
        +
        O(\Lambda^{-1}).
        \label{eq:z_middle_h}
\end{equation}
and insert it into \eqref{eq:geq}. The terms of order \(\sqrt{\Lambda}\) give
\begin{equation}
        2kg_0'
        +
        \left(
        k'
        +
        3\frac{f''}{f'}k
        \right)g_0
        =
        0,
\end{equation}
so that
\begin{equation}
        g_0
        =
        B_{\rm mid}(f')^{-7/4}.
        \label{eq:g0_middle}
\end{equation}
At the next order one obtains the equation for the first correction,
\begin{equation}
        2k g_1'
        +
        \left(
        k'
        +
        3\frac{f''}{f'}k
        \right)g_1
        =
        g_0''
        +
        3\frac{f''}{f'}g_0'
        +
        \frac{V_\beta}{f'}g_0 .
        \label{eq:h_middle}
\end{equation}
The solution for $g_1$ is
\begin{equation}
        g_1
        =
        B_{\rm mid}(f')^{-7/4}
        \left\{
        \frac{7}{48}\frac{f''}{(f')^{3/2}}
        +
        \frac{23\beta}{12}\frac{f}{(f')^{5/2}}
        +
        \int_0^\eta  J_\beta(y)\,dy
        \right\}.
        \label{eq:h_solution}
\end{equation}
where
\begin{equation}
        J_\beta
        =
        \frac{11}{48}f'''(f')^{-3/2}
        -
        \frac{1}{8}f^2(f')^{-1/2}
        +
        \frac{\beta}{4}(f')^{1/2}  
        -
        \frac{13\beta}{6}(f')^{-3/2}
        +
        \frac{115\beta}{24}f f''(f')^{-7/2}.
        \label{eq:Jbeta}
\end{equation}

It is clear from \eqref{eq:g0_middle} and \eqref{eq:h_solution} that the solution \eqref{eq:z_middle_h} cannot be extended all the way to the wall due to the singularities in $g_0$ and $g_1$. The wall condition must therefore be imposed in a separate inner layer.

\subsection{The wall approximation}
In the wall layer, Brown defines the inner local variable \(\zeta=\Lambda^{1/3}\eta\) that is $O(1)$ in this region. Hence, $\eta$ is small and $f$ can be replaced by its 
series expansion about $\eta=0$ so that \cite{coppel}
\begin{equation}
       f=\frac{a}{2}\eta^2-\frac{\beta}{6}\eta^3+O(\eta^5) ,\qquad f'=a\eta-\frac{\beta}{2}\eta^2+O(\eta^4).
        \label{smalletafexpansion}
\end{equation}
Brown then proceeds to write \eqref{eq:zEq} in terms of $\zeta$ and to approximate $A_\beta$, leaving to first order a Bessel equation. For \(\beta=0\), this leading Bessel solution is sufficient at the order needed for Brown's first correction. For $\beta>0$, however, the $\beta$-dependent terms of the wall expansion produce a relative $O(\Lambda^{-1/3})$ perturbation of the Bessel equation. This perturbation changes both the large-\(x\) phase and the large-\(x\) amplitude of the regular wall solution.

In order to address the wall layer in the present case, it is convenient to work with the middle phase variable
\begin{equation}
        x=\Lambda^{1/2}t
        =\frac23a^{1/2}\zeta^{3/2}
        \left\{1-\frac{3\beta}{20a}\Lambda^{-1/3}\zeta+\cdots\right\}.
        \label{xeta}
\end{equation}
It is also useful to set
\begin{equation}
        w=(f')^{7/4}z .
        \label{wtransf}
        \end{equation}
The transformations \eqref{eq:tauPhaseNorm} and \eqref{wtransf} reduce \eqref{eq:zEq} to normal
form exactly: with $x=\Lambda^{1/2}t$, the equation reads $\frac{d^2w}{dx^2}+Q\,w=0$, where
\begin{equation}
Q=\frac1\Lambda\left[\frac{A_\beta}{f'}-\frac74\frac{f'''}{(f')^2}
-\frac{7}{16}\frac{(f'')^2}{(f')^3}\right]
=1+\frac1\Lambda\left[\frac34\frac{f'''}{(f')^2}-\frac{7}{16}\frac{(f'')^2}{(f')^3}
-\frac{f^2}{4f'}+\frac{\beta}{2}\left(1-\frac{1}{(f')^2}\right)\right].
\label{eq:normalform}
\end{equation}
Substituting the wall series \eqref{smalletafexpansion} into \eqref{eq:normalform} gives 
\begin{equation}
Q=1-\frac1\Lambda\left[\frac{7}{16a\eta^{3}}+\frac{33\beta}{32a^{2}\eta^{2}}+\frac{33\beta^{2}}{32a^{3}\eta}+O(1)\right],
\label{eq:laurent}
\end{equation}
whose successive terms are of relative order $1$, $\Lambda^{-1/3}$ and $\Lambda^{-2/3}$ when
$\eta=O(\Lambda^{-1/3})$. The $\eta^{-3}$ term is $\beta$-independent and yields Brown's Bessel
operator; the $\eta^{-2}$ term, which carries the new correction, collects three $O(\beta)$ contributions --- from $f'''(0)=-\beta$, from the subleading terms of $f''$ and $f'$ in \eqref{smalletafexpansion}, and from the explicit $\beta$-term of $A_\beta$ --- adding up to $-33\beta/32a^{2}$.
Inverting \eqref{xeta} and substituting into \eqref{eq:laurent} makes it possible to write the normal-form equation in the wall layer as 
\begin{equation}
        w''+\left(1-\frac{7}{36x^2}
        -c_\beta\Lambda^{-1/3}x^{-4/3}
        -d_\beta\Lambda^{-2/3}x^{-2/3}+\cdots\right)w=0,
        \label{eq:wallEq}
\end{equation}
where primes now denote differentiation with respect to \(x\), and
\begin{equation}
        c_\beta=\frac{2^{1/3}3^{2/3}}{5}\frac{\beta}{a^{4/3}},
        \qquad
        d_\beta=\tfrac94\,c_\beta^{2}=\frac{2^{2/3}3^{10/3}}{100}\frac{\beta^2}{a^{8/3}}.
        \label{cdbeta}
\end{equation}

The regular solution is expanded as
\begin{equation}
        w=w_0+\Lambda^{-1/3}w_1+\cdots .
        \label{inner}
\end{equation}
Substituting into  \eqref{eq:wallEq} and equating powers of $\Lambda$ gives 
\begin{equation}
        \mathcal L_0 w_0=0,
        \label{Bessel01}
\end{equation}
\begin{equation}
        \mathcal L_0 w_1=c_\beta x^{-4/3}w_0,
           \label{Bessel1}     
\end{equation}
where 
\begin{equation}
\mathcal L_0 =d^2/dx^2+1-\frac{7}{36x^2} .
        \label{eq:Besselop}
\end{equation}
With
\begin{equation}
        \phi_1=x^{1/2}J_{2/3}(x),
        \qquad
        \phi_2=x^{1/2}Y_{2/3}(x),
\end{equation}
regularity gives \(w_0=c_1\phi_1\), and variation of constants gives
\begin{equation}
        w_1=\frac{\pi c_\beta c_1}{2}
        \left\{
        \phi_2\int_0^x y^{-4/3}\phi_1^2\,dy
        -\phi_1\int_0^x y^{-4/3}\phi_1\phi_2\,dy
        \right\}.
        \label{eq:w1}
\end{equation}

\section{Matching and eigenvalues}
\label{sec4}

We now need to match the solution of the middle
layer to the solutions in the outer and inner layers. We start with the solution in the outer layer, and note from \eqref{eq:outerosc} and \eqref{zetaout} that, for small $\xi$, 
\begin{equation}
 z_{\rm out}^{\xi\to0}(\xi)\sim B_{\rm out}\cos \Theta_{\xi\downarrow 0},
        \qquad 0<\xi<1,
        \label{Outerlower}
\end{equation}
where
\begin{equation}
        \Theta_{\xi\downarrow 0}
        =
        \frac{\pi}{2}\Lambda
        +\kappa_\beta\sqrt{\Lambda}
        -\sqrt{\Lambda}\eta
        -\frac{\pi}{4}
        +\frac{(\eta-\kappa_\beta)^3}{24\sqrt{\Lambda}}
        +\cdots .
        \label{eq:Xouter}
\end{equation}
This needs to be matched to the middle solution as $\eta\to\infty$. For large $\eta$ we have 
\begin{equation}
          Y(\eta)=-\sqrt{\Lambda} (\eta-\kappa_{1,\beta}+o(1))+\gamma ,
        \label{eq:Ylargeeta}
\end{equation}
where
\begin{equation}
        \kappa_{1,\beta}  =  \int_0^\infty \{1-f'^{1/2}\}\,d\eta 
\end{equation}
and
\begin{equation}
        \frac{g_1}{g_0}
        =
        -\frac{1}{24}(\eta-\kappa_\beta)^3
        +
        \delta_\beta
        +
        o(1),
        \qquad
        \eta\to\infty ,
        \label{eq:h_upper}
\end{equation}
where
\begin{equation}
        \delta_\beta
        =
        \lim_{\eta\to\infty}
        \Bigg[
        \frac{7}{48}\frac{f''(\eta)}{f'(\eta)^{3/2}}
        +
        \frac{23\beta}{12}
        \frac{f(\eta)}{f'(\eta)^{5/2}}  
        +
        \int_0^\eta J_\beta(y)\,dy
        +
        \frac{1}{24}(\eta-\kappa_\beta)^3
        \Bigg].
        \label{eq:delta_beta_explicit_h}
\end{equation}

Using \eqref{eq:Ylargeeta}, \eqref{eq:h_upper} and the asymptotic expansion of $f$ and its derivatives in \eqref{eq:z_middle_h} and writing the sum in  phase-amplitude form we find that the solution at the outer edge of the middle layer may be written
\begin{equation}
        z_{\rm mid}
        =
        B_{\rm mid}\sqrt{1+O(\Lambda^{-1})}\cos\left(\gamma-\Lambda^{1/2} (\eta-\kappa_{1,\beta})+\left((\eta-\kappa_\beta)^3/24
        -\delta_\beta
        \right)/\sqrt\Lambda
        \right)+
        O(\Lambda^{-1}).
        \label{eq:z_middle_outer}
\end{equation}
Comparison of \eqref{eq:z_middle_outer} and \eqref{Outerlower} yields the matching conditions
\begin{equation}
        \gamma
        =
        \frac{\pi}{2}\Lambda
        +
        \Delta_\beta\Lambda^{1/2}
        -
        \frac{\pi}{4}
        +
        \delta_\beta\Lambda^{-1/2}
        +
        s_1\pi
        +
        \cdots ,
        \label{eq:gamma_upper_h}
\end{equation}
and
\begin{equation}
        B_{\rm mid}=(-1)^{s_1}B_{\rm out}.
        \label{eq:Bmatch_h}
\end{equation}
where $s_1$ is an integer and
\begin{equation}
        \Delta_\beta
        =
        \kappa_\beta-\kappa_{1,\beta}
        =
        \int_0^\infty \{(f')^{1/2}-f'\}\,d\eta .
        \label{eq:Deltabeta_h}
\end{equation}

Matching the inner solution to the middle solution requires examining \eqref{inner} as $x\to\infty$ and \eqref{eq:z_middle_h} as $\eta\to 0$. Using the asymptotic expansion of the Bessel functions, the solution in the outer edge of the inner layer is, to the retained order,
\begin{equation}
\begin{split}
        w_0+\Lambda^{-1/3}w_1
        \sim{}&c_1\left(\frac{2}{\pi}\right)^{1/2}
        (1+\epsilon_\beta\Lambda^{-1/3}) \\
        &\times
        \cos\left(
        x-\frac{7\pi}{12}+\frac{7}{72x}
        +\frac{3c_\beta}{2}\Lambda^{-1/3}x^{-1/3}
        -\alpha_\beta\Lambda^{-1/3}
        \right),
\end{split}
\label{eq:wallLarge}
\end{equation}
where 
\begin{equation}
        \alpha_\beta=\frac{\pi c_\beta}{2}
        \int_0^\infty x^{-1/3}J_{2/3}(x)^2\,dx
        \label{eq:alphae}
\end{equation}
and
\begin{equation}
        \epsilon_\beta=-\frac{\pi c_\beta}{2}
        \int_0^\infty x^{-1/3}J_{2/3}(x)Y_{2/3}(x)\,dx.
        \label{eq:epsilon}
\end{equation}

Using \eqref{xeta}, the phase in \eqref{eq:wallLarge} can be written in terms of $\eta$ as 
\begin{equation}
%
x-\frac{7\pi}{12}-\alpha_\beta\Lambda^{-1/3}+\frac{7}{48\sqrt a}
        \Lambda^{-1/2}\eta^{-3/2}+\frac{59\beta}{64a^{3/2}}
        \Lambda^{-1/2}\eta^{-1/2}+O(\Lambda^{-1/2}\eta^{1/2}).
        \label{eq:largeinnerphase}
\end{equation}

Using the series expansion of $f$ and its derivatives near the origin, we obtain
$J_\beta(\eta)=O(\eta^{-1/2})$ as $\eta\downarrow 0$ and, thus, $
\int_0^\eta J_\beta(y)\,dy = O(\eta^{1/2})$,
so that the integral term in $g_1$ tends to zero and does not contribute to the
singular lower-overlap phase. The middle solution therefore behaves as
\begin{equation}
        z_{\rm mid}(\eta)
        \sim
        g_0(\eta)
        \cos\left[
        x-\gamma
        +
        \Lambda^{-1/2}\left(
\frac{7}{48\sqrt a}\eta^{-3/2}
        +
        \frac{59\beta}{64a^{3/2}}\eta^{-1/2}
        +
        O(\eta^{1/2})\right)        
        \right].
        \label{eq:middle_lower_phase_h}
\end{equation}
Comparison with  \eqref{eq:largeinnerphase} yields the matching condition
\begin{equation}
        \gamma = \frac{7\pi}{12} + \alpha_\beta\Lambda^{-1/3} + s_2\pi + O(\Lambda^{-2/3}),
        \label{eq:gamma_lower}
\end{equation}
where \(s_2\) is another integer. The remainder in \eqref{eq:gamma_lower} is $O(\Lambda^{-2/3})$ because the inner expansion \eqref{eq:wallEq} proceeds in powers of $\Lambda^{-1/3}$, so no $\Lambda^{-1/2}$ constant arises on the wall side. Combining \eqref{eq:gamma_upper_h} and \eqref{eq:gamma_lower} we have 
\begin{equation}
        \frac{\pi}{2}\Lambda
        +\Delta_\beta\Lambda^{1/2}
        -\frac{5\pi}{6}
        -\alpha_\beta\Lambda^{-1/3}
        +\delta_\beta\Lambda^{-1/2}
        +\cdots=s\pi ,
        \label{eq:phaseCondition}
\end{equation}
where $s=s_2-s_1$. Solving for \(\lambda=\Lambda+1/2+2\beta\), and writing \(s=n-1\), where $n=1,2,3,\ldots$ is the eigenmode index, gives
\begin{equation}
\begin{split}
        \lambda_n(\beta)=
        2s-\frac{2\Delta_\beta}{\pi}(2s)^{1/2}
        +\frac{13}{6}+2\beta+\frac{2\Delta_\beta^2}{\pi^2}
        +E_\beta(2s)^{-1/3}+F_\beta(2s)^{-1/2}
        +O(s^{-2/3}),
\end{split}
\label{eq:lambdaGeneral}
\end{equation}
where
\begin{equation}
        E_\beta=\frac{2\alpha_\beta}{\pi},
        \qquad
        F_\beta=-\frac{2\delta_\beta}{\pi}
        -\frac{5\Delta_\beta}{3\pi}
        -\frac{\Delta_\beta^3}{\pi^3}.
\end{equation}
Using \eqref{cdbeta} and standard Bessel integral identities \cite{Watson1966},
$\alpha_\beta$ can be written as
\begin{equation}
        \alpha_\beta=\frac{3^{5/3}\pi}{10\Gamma(2/3)^2}\frac{\beta}{a^{4/3}}.
        \label{eq:alpha2}
\end{equation}

For \(\beta=0\), \(E_\beta=0\) and Brown's form is recovered. For \(\beta=1/2\), the numerical constants are
\begin{equation}
        a=0.92768,
        \qquad
        \Delta_{1/2}=0.29493,
        \qquad
        \delta_{1/2}=0.15525 ,
        \qquad
        \alpha_{1/2}= 0.59085.
\end{equation}

Thus
\begin{equation}
        \lambda_n(1/2)
        \simeq
        2s-0.187756(2s)^{1/2}+3.18429
        +0.37615(2s)^{-1/3}-0.25613(2s)^{-1/2}.
        \label{eigenveq}
\end{equation}

To assess \eqref{eigenveq}, we computed reference eigenvalues by solving the
Chen-Libby problem directly. The Falkner-Skan profile is obtained by
shooting~\eqref{eq:FS} for $a=f''(0)$ so that $f'(\infty)=1$. The reduced equation~\eqref{eq:zEq} is then solved by two-sided shooting. The wall-regular solution, fixed by $z(0)=1$, $z'(0)=\beta/a$ from~\eqref{eq:zWallBC} and started
from $\eta=10^{-4}$ through the regular series to bypass the $f'\to0$ wall singularity, is
integrated outward to a match point $\eta_m=\eta_0-\delta_m$ in the oscillatory interior just
below the turning point $\eta_0=\kappa_\beta+2\Lambda^{1/2}$. The exponentially small solution
is launched from a far-field cut-off $\eta_\infty=\eta_0+\delta_\infty$ with the decaying WKB
slope $z'/z=-\sqrt{-A_\beta(\eta_\infty)}$ and integrated inward---the numerically stable
direction for the subdominant branch---to $\eta_m$. The eigenvalue is the zero of the Wronskian
$z_Lz_R'-z_L'z_R$ at $\eta_m$, bracketed about the estimate~\eqref{eigenveq} and refined
to $10^{-10}$. All integrations were carried out in IEEE double precision with an adaptive-step
solver at relative/absolute tolerances of $10^{-11}/10^{-13}$ (the Falkner--Skan profile to
$3\times10^{-13}$); no extended-precision arithmetic was needed for the six--seven significant
figures of Tables~\ref{tab:eig_beta_half} and~\ref{tab:thirty_modes}. The eigenvalues are independent of these
choices to the reported figures: varying the match offset over $1\le\delta_m\le3.5$ and the
far-field offset over $3\le\delta_\infty\le9$ leaves $\lambda_n$ unchanged, so neither the
domain truncation nor the match placement limits the reference accuracy---a prerequisite for
using it to gauge the earlier tabulations at high $n$.

The results are collected in Table~\ref{tab:eig_beta_half}. The numerical eigenvalues coincide
with those tabulated by Chen and Libby~\cite{ChenLibby1968} to within their three-digit
precision, and the asymptotic formula~\eqref{eigenveq} reproduces both, the relative error
falling from $0.34\%$ at $n=5$ to $0.03\%$ at $n=20$. The new term of $O(s^{-1/3})$ is significant where it is tested: at $n=5$ it contributes $E_\beta(2s)^{-1/3}=0.188$, a $1.8\%$ shift, so omitting it turns the $+0.34\%$ agreement into
a $-1.4\%$ error; even at $n=20$ it still contributes $0.112$ ($0.28\%$ of $\lambda_{20}$),
whose omission would leave a $-0.25\%$ error where \eqref{eigenveq} achieves $+0.03\%$. The errors of \eqref{eq:lambdaGeneral} themselves decay at the predicted rate: the absolute error falls from $0.036$ at $n=5$ to $0.0094$ at $n=30$, a ratio $3.8$ consistent with the $(2s)^{-2/3}$ remainder $\big((58/8)^{2/3}=3.75\big)$.

\begin{table}[H]
  \centering
  \begin{tabular}{c c c c c c}
    \toprule
    $n$ & $s=n-1$ & $\lambda_n$ (num.) & $\lambda_n$~\eqref{eigenveq}
       & error (\%) & $\lambda_n$~\cite{ChenLibby1968} \\
    \midrule
     5 &  4 & 10.7148 & 10.7508 & 0.336\% & 10.715 \\
    10 &  9 & 20.4498 & 20.4709 & 0.103\% & 20.450 \\
    15 & 14 & 30.2506 & 30.2663 & 0.052\% & 30.251 \\
    20 & 19 & 40.0846 & 40.0972 & 0.031\% & 40.088 \\
    25 & 24 & 49.9393 & 49.9500 & 0.021\% & ---    \\
    30 & 29 & 59.8085 & 59.8179 & 0.016\% & ---    \\
    \bottomrule
  \end{tabular}
  \caption{Eigenvalues for $\beta=1/2$. ``num.'' is a numerical
    reference computed by two-sided shooting;
    $\lambda_n$~\eqref{eigenveq} is the asymptotic formula; the last column
    lists the values tabulated by Chen and Libby~\cite{ChenLibby1968} (available
    up to $n=20$). }
  \label{tab:eig_beta_half}
\end{table}

\section{Norms}
\label{sec5}

The norms of the Chen-Libby eigenfunctions are defined as \cite{ChenLibby1968} 
\begin{equation}
        C_n
        =
        \int_0^\infty (f')^4 e^F \left(({N_n}/{f'})'\right)^2\,d\eta ,
\end{equation}
where $F=\int_0^\eta f(\mu)d\mu$. Tracing through the transformations \eqref{ztransf},  \eqref{eq:tauPhaseNorm} and \eqref{wtransf}, this can be written as 
\begin{equation}
        C_n
        =
        \int_0^\infty (f')^4 z_n(\eta)^2\,d\eta = \int_0^\infty w_n(t)^2\dd t .
        \label{eq:normzBrown}        
\end{equation}

In order to compute the norm we need to start by fixing the amplitude constants $c_1$, $B_{\rm mid}$ and $B_{\rm out}$. The last two are related by \eqref{eq:Bmatch_h}, while the first one is fixed by the wall boundary behavior of the Chen-Libby eigenfunctions. Assuming that \(N_n''(0)=1\), we have \(N_n\sim \eta^2/2\) and \(f'\sim a\eta\), from where we get
\begin{equation}
        z_n(\eta)\sim \frac{1}{2a},
        \qquad
        w_n(t)\sim S_\beta t^{7/6},
        \qquad
        S_\beta=\frac12\left(\frac32\right)^{7/6}a^{1/6}.
        \label{esebeta}
\end{equation}
This fixes the constant in the wall solution:
\begin{equation}
        c_1=2^{2/3}\Gamma\left(\frac53\right)S_\beta\Lambda^{-7/12}.
        \label{c1}
\end{equation}
The inner solution at the overlap with the middle layer \eqref{eq:wallLarge} has amplitude $c_1\sqrt{\frac2\pi}(1+\epsilon_\beta\Lambda^{-1/3})$. Given \eqref{eq:middle_lower_phase_h} and that $w_{\rm mid}=f'^{7/4}z_{\rm mid}$, the amplitude of the middle solution is, thus, 
\begin{equation}
        B_{\rm mid}=A_0S_\beta\Lambda^{-7/12}
        \left(1+\epsilon_\beta\Lambda^{-1/3}+O(\Lambda^{-2/3})\right),
        \qquad
        A_0=2^{2/3}\Gamma\left(\frac53\right)\left(\frac{2}{\pi}\right)^{1/2}.
        \label{BmidA0}
\end{equation}
This is equal, up to a sign, to the amplitude constant of the outer solution \eqref{eq:outerosc}. The sign is irrelevant as far as the computation of the norm is concerned, so henceforth we drop the labels ${\rm mid/out}$ and refer to the absolute value of the amplitude constant of the middle and outer solutions as $B$ with a subscript $n$ to distinguish the particular eigenfunction. 

The integral         \eqref{eq:normzBrown} should be a priori divided into four patches: the inner (wall) and middle layers, the oscillatory part of the outer layer up to the turning point, and the exponentially decaying tail beyond it. The exact norm would require the appropriate eigenfunction in each patch. We instead discard the contribution of the region beyond the turning point and replace the wall/middle integrand by the outer amplitude \eqref{eq:outerosc} continued to $t=0$. In this approximation, the wall information is incorporated through the amplitude constant $B_n$ of \eqref{BmidA0}, fixed by matching across the wall and middle layers. 

We thus compute the norm by the following integral
\begin{equation}
C_n\simeq\int_0^{t_+}\langle w_n^2\rangle\,dt\simeq \frac{B_n^2}{2}\int_0^{t_+}\!\mathcal E(t)\,dt
\label{Cref}
\end{equation}
where $\langle\,\cdot\,\rangle$ denotes the average over the fast oscillation of the phase in \eqref{eq:outerosc}, 
\begin{equation}
\langle w_n^2\rangle \simeq \tfrac12 B_n^2\,\mathcal E(t),\qquad
\mathcal E(t)=\Big[1-\frac{(t-\Delta_\beta)^2}{4\Lambda_n}\Big]^{-1/2}
\end{equation}
and
\begin{equation}
t_+=\Delta_\beta+2\Lambda_n^{1/2}
\label{tmas}
\end{equation}
is the image of the turning point $\eta_0$ under \eqref{eq:tauPhaseNorm} and the point at which $\mathcal E$ is singular. The integral \eqref{Cref} is elementary, and gives
\begin{equation}
C_n
\simeq\frac{B_n^2}{2}\Big(\pi\Lambda_n^{1/2}+\Delta_\beta+\frac{\Delta_\beta^3}{24\Lambda_n}+\cdots\Big).
\label{eq:Cref2}
\end{equation}

The choice of \eqref{eq:outerosc} as the integrand deserves a word. None of the three layer solutions is uniformly valid, so we continue one of them outside its own layer, and the outer oscillatory form is the natural candidate: it alone covers the bulk of the range of integration. In the variable $t$, the turning point lies at $t_+\simeq2\Lambda_n^{1/2}$, while the middle layer occupies $t=O(1)$ and the wall layer only $t=\Lambda_n^{-1/2}x=O(\Lambda_n^{-1/2})$. The cost of the continuation is then easy to weigh. Both the true integrand and its replacement are $O(B_n^2)$ throughout the displaced regions: in the middle layer, $\mathcal E(t)=1+O(\Lambda_n^{-1})$ and $w_{\rm mid}=B_n\cos Y+O(\Lambda_n^{-1/2})$, so the amplitude is $B_n$ to leading order, while in the wall layer, where the true amplitude does depart from $B_n$, both remain bounded and the interval is only $O(\Lambda_n^{-1/2})$ long.
Their difference, integrated, is therefore at most $O(B_n^2\Lambda_n^{-1/2})$, to be compared with the leading $O(B_n^2\Lambda_n^{1/2})$ of \eqref{eq:Cref2}: a relative error $O(\Lambda_n^{-1})$, beyond the orders retained. Integrating the Bessel or middle solutions instead would resolve regions that carry a vanishing fraction of the norm.  Equation~\eqref{Cref} is, in this sense, an approximation in the spirit of the norm estimate of Kemp~\cite{Kemp1970}.

Two approximations remain that are not obviously innocuous, both at the outer end: the truncation of the integral at $t_+$, and the use of the envelope $\mathcal E$ up to $t_+$, where it is singular. We do not carry out a careful check of the effect of these two outer-end approximations on the retained orders, so the norm expansions presented below in \eqref{eq:normLeading} and \eqref{eq:normSecond} are, at best, semi-asymptotic rather than rigorously derived and will be indirectly validated by comparison with numerical computations, which support the idea that the neglected contributions affect only orders beyond those retained in the final norm estimates.

We take \eqref{eq:Cref2} as the representative of the norm to the orders retained, and assess the accuracy of this step a~posteriori against the numerical norms of Table~\ref{tab:norms_beta_half}. Substituting \eqref{BmidA0} into \eqref{eq:Cref2} yields
\begin{equation}
        C_n(\beta)
        \sim
        K_\beta\Lambda_n^{-2/3}
        \left(
        1
        +
        P_\beta\Lambda_n^{-1/3}
        +
        \frac{\Delta_\beta}{\pi}\Lambda_n^{-1/2}
        \right),
        \label{eq:normLeading}
\end{equation}
where
\begin{equation}
        K_\beta
        =
        S_\beta^2\,2^{4/3}\Gamma\left(\frac53\right)^2,
        \qquad
        P_\beta=2\epsilon_\beta .
\end{equation}
Using \eqref{cdbeta} and standard Bessel integral identities \cite{Watson1966},
we get
\begin{equation}
        \epsilon_\beta=\frac{\pi }{2}
        \frac{3^{7/6}\beta}{5\Gamma(2/3)^2a^{4/3}} .
        \label{eq:epsilon2}
\end{equation}
For \(\beta=1/2\),
\begin{equation}
K_{1/2}=1.28960,\qquad P_{1/2}=0.68226,\qquad
\frac{\Delta_{1/2}}{\pi}=0.09388, 
\end{equation}
which yields a first estimate
\begin{equation}
        C_n(1/2)
        \sim
        1.28960\,\Lambda_n^{-2/3}
        \left\{
        1
        +0.68226\,\Lambda_n^{-1/3}
        +0.09388\,\Lambda_n^{-1/2}
        \right\}.
        \label{eq:normbeta12first}
\end{equation}
Measured against Chen and Libby's tabulated norms \cite{ChenLibby1968}, the leading estimate \eqref{eq:normbeta12first}  is off by around 8\% at the 20th eigenmode. As Table \ref{tab:norms_beta_half} shows, however, much of that gap lies in the early tabulation rather than in the estimate: against our own numerical computations, \eqref{eq:normbeta12first} is low by only about 4\% there. 

We can improve the agreement by computing the next term, of $O(\Lambda_n^{-2/3})$. This requires
the second-order inner correction $w_2$ of the expansion \eqref{inner}, constructed in
Appendix~\ref{appA}. The second-order source generates two projection integrals of different
character: a convergent one, which defines the constant $\epsilon_{2\beta}$ of
\eqref{eq:eps2}, and a slowly growing one, which modulates the phase of the inner oscillation
rather than its amplitude. Only the former enters the envelope, and the large-$x$ amplitude of
the inner solution \eqref{eq:wallLarge} is promoted to
\begin{equation}
c_1\sqrt{\frac2\pi}\,\left(1+\epsilon_\beta\Lambda_n^{-1/3}
+\Big(\epsilon_{2\beta}+\tfrac12\alpha_\beta^2\Big)\Lambda_n^{-2/3}+O(\Lambda_n^{-1})\right),
\label{2ndorderamplitude}
\end{equation}
where $\tfrac12\alpha_\beta^2$ is the second-order contribution of the first-order sine admixture $\alpha_\beta$ of \eqref{eq:alphae} and $\epsilon_{2\beta}$ is the second-order cosine amplitude of the inner solution defined by~\eqref{eq:eps2} (see Appendix~\ref{appA}). 
Eq. \eqref{2ndorderamplitude} yields, through the appropriate modification of the amplitude constant $B_n$, the following improved norm estimate 
\begin{equation}
        C_n(\beta)
        \sim
        K_\beta\Lambda_n^{-2/3}
        \left(1+P_\beta\Lambda_n^{-1/3}
        +\frac{\Delta_\beta}{\pi}\Lambda_n^{-1/2}
        +R_\beta\Lambda_n^{-2/3}
        \right),
        \label{eq:normSecond}
\end{equation}
where
\begin{equation}
R_\beta=\epsilon_\beta^2+2\Big(\epsilon_{2\beta}+\tfrac12\alpha_\beta^2\Big).
        \label{eq:qbigbeta}
\end{equation}

For \(\beta=1/2\), 
\begin{equation}
\epsilon_{2\beta}=4.6\times10^{-4},\qquad R_{1/2}=0.46639 .
\end{equation}

Hence,
\begin{equation}
        C_n(1/2)
        \sim
        1.28960\,\Lambda_n^{-2/3}
        \left\{
        1
        +0.68226\,\Lambda_n^{-1/3}
        +0.09388\,\Lambda_n^{-1/2}
        +0.46639\,\Lambda_n^{-2/3}
        \right\}.
        \label{eq:normbeta12}
\end{equation}

The norm estimate \eqref{eq:normbeta12} still differs significantly from Chen-Libby data. To assess whether this is a problem of the approach or of the reference data, we computed the norms for the numerical eigenfunctions obtained in Section \ref{sec4}. The norms
$C_n=\int_0^\infty (f')^4 z_n^2\,\mathrm{d}\eta$ are evaluated by quadrature after rescaling the eigenfunction to $z_n(0)=1/(2a)$, i.e. $N_n''(0)=1$. The composite eigenfunction---the wall-regular and decaying branches matched at $\eta_m$---is
used throughout, and the quadrature is carried out to $\eta_\infty$, where $(f')^4z_n^2$ has
fallen below the tolerance. The results are converged to better than $10^{-5}$. (Appendix \ref{appB} compiles the eigenvalues and norms for the first 30 Chen-Libby eigenfunctions for $\beta=1/2$ computed with the two-sided shooting method). 

Table~\ref{tab:norms_beta_half} compares the  estimates~\eqref{eq:normbeta12first}
and~\eqref{eq:normbeta12} with the numerical reference results. The leading
estimate~\eqref{eq:normbeta12first} lies below the reference --- by about 9\% at $n=5$, decreasing monotonically to roughly 3\% by $n = 30$ --- as expected of an asymptotic series whose first neglected term is $O(\Lambda_n^{-2/3})$. Including the
$O(\Lambda_n^{-2/3})$ wall-amplitude term~\eqref{eq:normbeta12} roughly halves the
error, bringing it to about $2\%$ at $n=5$ and below $0.7\%$ by $n=20$. The
norms tabulated by Chen and Libby~\cite{ChenLibby1968} agree with the present
numerical reference to about $0.2\%$ at $n=5$ but lie progressively above it as
$n$ increases, the difference reaching about $4\%$ at $n=20$; this growth with
mode number is consistent with reduced quadrature accuracy of the early
tabulation for the more rapidly oscillating high-order eigenfunctions. The slow
convergence of the norm series itself reflects that its expansion parameter
$\Lambda_n^{-1/3}$ is only moderately small over the tabulated modes
($\Lambda_n^{-1/3}\approx0.49$ at $n=5$, $0.29$ at $n=20$). That including the $O(\Lambda_n^{-2/3})$ term \eqref{eq:normbeta12} reduces the error several-fold and leaves a remainder decaying faster still is the behavior of a genuine asymptotic expansion. What is more, the error of \eqref{eq:normbeta12} in Table~\ref{tab:norms_beta_half} can be seen to exhibit a decay rate $\approx\Lambda_n^{-0.87}$, which is compatible with the cross term $O(\Lambda_n^{-5/6})$ arising on expanding
$(1+P_\beta\Lambda_n^{-1/3}+R_\beta\Lambda_n^{-2/3})\big(1+(\Delta_\beta/\pi)\Lambda_n^{-1/2}\big)$ --- a rate $O(\Lambda_n^{-1/3})$ or $O(\Lambda_n^{-2/3})$ terms could not produce, so the numerics independently place the residue below $\Lambda_n^{-2/3}$ and confirm a~posteriori that the approximation~\eqref{Cref} does
not disturb the retained orders. The same $\approx\Lambda_n^{-5/6}$ decay also certifies that the slowly growing projection
integral of Appendix~\ref{appA} does not enter the norm at the retained order, since an
unaccounted $O(\Lambda_n^{-2/3})$ term in $R_\beta$ would produce an $O(\Lambda_n^{-2/3})$
residue rather than the observed faster decay.

\begin{table}[H]
  \centering
  \begin{tabular}{c c c c c c c c}
    \toprule
    $n$ & $\lambda_n$ & $C_n$ (num.) & first~\eqref{eq:normbeta12first} & err.\ (\%) & second~\eqref{eq:normbeta12} & err.\ (\%) & $C_n$~\cite{ChenLibby1968} \\
    \midrule
         5 & 10.7148 & 0.43927 & 0.39796 & $-9.40$ & 0.42909 & $-2.32$ & 0.440 \\
        10 & 20.4498 & 0.24680 & 0.23178 & $-6.09$ & 0.24368 & $-1.26$ & 0.248 \\
        15 & 30.2506 & 0.17882 & 0.17042 & $-4.70$ & 0.17725 & $-0.88$ & 0.182 \\
        20 & 40.0846 & 0.14303 & 0.13745 & $-3.90$ & 0.14206 & $-0.68$ & 0.149 \\
        25 & 49.9393 & 0.12059 & 0.11652 & $-3.38$ & 0.11993 & $-0.55$ & ---   \\
        30 & 59.8085 & 0.10506 & 0.10191 & $-3.00$ & 0.10457 & $-0.47$ & ---   \\
    \bottomrule
  \end{tabular}
  \caption{Norms for $\beta=1/2$. ``num.'' is the numerical
    reference (two-sided shooting); ``first'' and ``second'' are the estimates~\eqref{eq:normbeta12first}
    and~\eqref{eq:normbeta12}, with signed relative errors against the reference.
    The estimates are evaluated at the numerical eigenvalues of
    Table~\ref{tab:eig_beta_half}. The last column lists the values tabulated by
    Chen and Libby~\cite{ChenLibby1968}, which agree with the numerical reference
    to about $0.2\%$ at $n=5$ but lie progressively above it at larger $n$,
    reaching about $4\%$ at $n=20$.}
  \label{tab:norms_beta_half}
\end{table}

For the Blasius case, \(\beta=0\) and, thus, \(c_0=d_0=P_0=0\). The Blasius norm in Libby and Fox's convention is obtained by dividing \eqref{eq:normLeading} by \(a_0=f_0''(0)\). 
Thus,
\begin{equation}
C^{\rm LF}_n(0)\sim 2.18862\,\Lambda_n^{-2/3}\left(1+0.13526\,\Lambda_n^{-1/2}\right),
\qquad \Lambda_n=\lambda_n-\tfrac12 .        
        \label{normBlasius}
\end{equation}
Two-sided shooting was applied in this case as well in order to generate reference data. The wall condition~\eqref{eq:zWallBC} reduces to $z'(0)=0$, but the rest of the procedure is the same. The results are compiled in Table   \ref{tab:norms_blasius}. Norms are reported in the Libby-Fox convention $C^{LF}_n=C_n/a_0$, $a_0=f''_0(0)=0.46960$, for comparison
with Libby~\cite{Libby1965}. The numerically obtained eigenvalues reproduce
Brown's asymptotic values~\cite{Brown1968} to within $10^{-3}$ over the shown range, giving an independent check of the $\beta=0$ eigenvalue expansion.

\begin{table}[H]
  \centering
  \begin{tabular}{r r r r r r}
    \toprule
    $n$ & $\lambda_n$ & $C^{LF}_n$ (num.) & Eq.~\eqref{normBlasius}
        & err.\ (\%) & $C^{LF}_n$~\cite{Libby1965} \\
    \midrule
     5 &  9.4144 & 0.54506 & 0.53213 & $-2.37$ & 0.5444 \\
    10 & 19.0397 & 0.32583 & 0.32226 & $-1.10$ & 0.330  \\
    15 & 28.7591 & 0.24363 & 0.24190 & $-0.71$ & 0.246  \\
    20 & 38.5247 & 0.19883 & 0.19780 & $-0.52$ & 0.197  \\
    25 & 48.3193 & 0.17007 & 0.16937 & $-0.41$ & ---    \\
    \bottomrule
  \end{tabular}
  \caption{Norms for the Blasius case ($\beta=0$), Libby-Fox convention
    $C^{LF}_n=C_n/a_0$. ``num.'' is the numerical reference
    (two-sided shooting); Eq.~\eqref{normBlasius} is the asymptotic estimate, with
    signed relative error against the reference. Since $\beta=0$ gives $P_0=R_0=0$,
    the first and second estimates coincide. The last column lists the values
    tabulated by Libby~\cite{Libby1965}.}
  \label{tab:norms_blasius}
\end{table}

As for $\beta=1/2$, the estimate~\eqref{normBlasius} underpredicts the numerical reference value by a few per cent and improves monotonically with $n$, falling below 0.5\% by $n=25$. The differences from Libby's tabulated
values~\cite{Libby1965} are larger and non-monotonic, reflecting the precision
of that early data rather than the expansion: the present numerics indicate that
the asymptotic estimate, not the 1965 tabulation, is the more accurate of the two
at the modes shown. As in the $\beta=1/2$ case, the errors of \eqref{normBlasius} in Table~\ref{tab:norms_blasius} decay at a rate
$\approx\Lambda_n^{-1}$, compatible with a first neglected term $O(\Lambda_n^{-1})$; this provides the same a~posteriori confidence to the two retained terms.

\section{Extension to adverse pressure gradients \texorpdfstring{$(\beta<0)$}{(beta<0)}}
\label{sec6}

For adverse-pressure Falkner--Skan profiles ($\beta<0$), the matched--asymptotic framework of Sections \ref{sec3}--\ref{sec4} can be formally extended to the upper branch. The changes relative to the favorable-gradient case are summarized below. Since no separate adverse-gradient numerical benchmark is included here, this section should be regarded as a formal extension of the asymptotic construction, rather than as an independently validated numerical result.

\paragraph{Upper branch.} The present method extends to the   \emph{upper} solution of the Falkner-Skan equation for which
        $f''(0)=a>0$ and $f'(\eta)>0$ for all $\eta>0$.  This branch exists for $\beta_0<\beta<0$, with $\beta_0\simeq -0.1988$, and coalesces with the lower branch as $\beta\to\beta_0^+$.

\begin{enumerate}
  \item \textbf{Wall layer.}  The expansion \eqref{smalletafexpansion} is unchanged, but now $\beta<0$ so the $c_\beta$ perturbation term in
        eq.~\eqref{eq:wallEq} enters with the \emph{opposite sign}. 
        Consequently $c_\beta\;(\text{eq.\;\eqref{cdbeta}})<0$, leading to a \emph{negative} $E_\beta$ and hence a
        \textbf{negative $s^{-1/3}$ shift} in the eigenvalue expansion
        \eqref{eq:lambdaGeneral}.

  \item \textbf{Middle and outer layers.}  Because
        $f'(\eta)$ stays positive, the turning-point structure and WKB
        matching are unaltered.  The constants
        $\kappa_\beta$, $\Delta_\beta$, $\alpha_\beta$, $\delta_\beta$
        must be re-evaluated with the new
        $f(\eta;\beta)$ but the derivations are identical.

  \item \textbf{Spectrum.}  On the upper branch the eigenvalues remain real, discrete and positive for all $n$ \cite{ChenLibby1968}, so disturbances decay monotonically downstream. The eigenvalues follow \eqref{eq:lambdaGeneral} with $E_\beta\propto c_\beta\propto\beta<0$: the leading wall correction flips
sign, while the $\Lambda^{-1/2}$ term and the overall hierarchy are unchanged.

        \item \textbf{Norms.} Equations~\eqref{eq:normLeading}-\eqref{eq:normSecond} carry over; $\alpha_\beta,\epsilon_\beta\propto\beta$ change sign, while
$\Delta_\beta,\delta_\beta$ are re-evaluated with the new profile.

  \item \textbf{Caveat.} 
The fixed-$\beta$ expansion is expected to remain valid away from separation; as
$\beta\to\beta_0^+$ one has $a\to0$, and the wall constants
$c_\beta,d_\beta\propto a^{-4/3},a^{-8/3}$ diverge, signaling the breakdown of the
fixed-$\beta$ ordering at the separation point.      
\end{enumerate}

\paragraph{Lower branch ($f''(0)<0$).}
On the lower branch $a<0$ and the base profile reverses near the wall:
$f'(\eta)<0$ for $0<\eta<\eta_c$ and $f'(\eta)>0$ for $\eta>\eta_c$, with $f'(\eta_c)=0$.
In the reverse-flow sublayer, $\sqrt{f'}$ is imaginary, so the wall variable $x$ becomes
imaginary, and the regular inner solution is a modified Bessel function ($I_\nu$, $K_\nu$)
rather than the oscillatory $J_\nu$, $Y_\nu$. This does not render the spectrum complex:
the zero of $f'$ at $\eta_c$ is an additional internal turning point at which the
non-oscillatory wall solution connects to the oscillatory middle region, so a matched
construction on the lower branch must pass through both $\eta_c$ and the outer turning
point. Chen and Libby~\cite{ChenLibby1968}, by a Sturm--Liouville argument that splits the
range at $\eta_c$, showed that the lower-branch spectrum is nevertheless real, comprising
both positive and negative eigenvalues separated by a $\beta$-dependent gap that widens
without bound as $\beta\to\beta_0^{+}$; the flow is correspondingly spatially unstable.
A matched-asymptotic treatment of this two-turning-point structure lies beyond the scope
of this paper.


\section{Conclusions}
\label{sec:conclusions}

We have extended Brown's matched-asymptotic construction for the Blasius problem
($\beta=0$) to Falkner-Skan profiles with positive wall shear $a>0$: the single
branch at favorable gradient $\beta>0$, and the upper branch of the adverse range
$-0.1988<\beta<0$.

\begin{enumerate}
  \item A new wall–layer perturbation proportional to
        $\beta\Lambda^{-1/3}x^{-4/3}$ modifies the inner Bessel problem and
        produces the \emph{first} correction $O(s^{-1/3})$ in the
        large–mode eigenvalue expansion~\eqref{eq:lambdaGeneral}.  Brown's
        $O(s^{-1/2})$ term survives but with a $\beta$–dependent constant
        obtained by overlap matching.
  \item Closed-form expressions were obtained for the matching constants:
$\alpha_\beta$ as an elementary Bessel integral, and $\Delta_\beta,\delta_\beta$ as
integrals over the Falkner-Skan profile; all reduce smoothly to Brown's values as
$\beta\to0$.
\item The same composite eigenfunction yields compact estimates of the Chen--Libby
norms, \eqref{eq:normLeading}--\eqref{eq:normSecond}. For $\beta=1/2$ the leading estimate \eqref{eq:normbeta12first} is accurate to a few per cent, and the second wall-amplitude term \eqref{eq:normbeta12} roughly halves the error, reaching
sub-percent accuracy by $n\simeq 20$. The remaining errors decay at the rate set by
the first neglected term of each estimate --- $\approx\Lambda_n^{-2/3}$ for \eqref{eq:normLeading} and
$\approx\Lambda_n^{-5/6}$ for \eqref{eq:normSecond} --- which is the behaviour of a genuine asymptotic
expansion and, in the absence of a formal error bound, is what justifies the norm
integral \eqref{Cref} a~posteriori. The numerical eigenvalues reproduce the classical
tabulations of Chen and Libby \cite{ChenLibby1968} and Libby \cite{Libby1965} to three-digit accuracy; the
corresponding tabulated norms agree at the lowest modes but lie a few per cent above
the present reference at the higher modes shown.

\end{enumerate}

The construction is limited to flows with $a>0$. On the lower branch, the base profile reverses near the wall, introducing an additional internal turning point; the spectrum is nonetheless real, with both positive and negative eigenvalues~\cite{ChenLibby1968}, and its asymptotic treatment would require a separate two-turning-point matched analysis.

\section*{Data Availability}
The data that support the findings of this study are available from the corresponding authors upon reasonable request.

\section*{Funding}
The research described in this paper has been supported by INTA under grant IDATEC (IGB21001).

\begin{appendices}

\section{The second-order wall amplitude}
\label{appA}

The $O(\Lambda_n^{-2/3})$ term of \eqref{eq:normSecond} requires the second-order inner
correction $w_2$ in the expansion \eqref{inner}, which obeys
\begin{equation}
\mathcal L_0 w_2=F_2,\qquad F_2=c_\beta x^{-4/3}w_1+d_\beta x^{-2/3}w_0 .
\label{Bessel2}
\end{equation}
The regular solution follows by variation of constants; with the Wronskian
$\mathcal W\{\phi_1,\phi_2\}=2/\pi$, we have
\begin{equation}
w_2 = \frac{\pi}{2}\left[\phi_2(x)\int_0^x \phi_1 F_2\,dy
-\phi_1(x)\int_0^x \phi_2 F_2\,dy\right].
\label{eq:w2sol}
\end{equation}
As $x\to\infty$, $\phi_1\sim\sqrt{2/\pi}\,\cos\Phi$ and $\phi_2\sim\sqrt{2/\pi}\,\sin\Phi$,
with $\Phi=x-\tfrac{7\pi}{12}+\cdots$. The two integrals in \eqref{eq:w2sol} behave
differently. Since $\phi_2F_2$ has no non-integrable non-oscillatory part at large $x$ (its leading term is
$\propto x^{-2/3}\phi_1\phi_2\sim x^{-2/3}\sin2\Phi$), the second integral converges to a well-defined constant 
\begin{equation}
  \int_0^x \phi_2 F_2\,dy \;\longrightarrow\; \int_0^\infty \phi_2 F_2\,dy
  \;=\; -\frac{2c_1}{\pi}\,\epsilon_{2\beta}.
\label{eq:eps2}
\end{equation}
where $\epsilon_{2\beta}$ is a new constant parameter. In contrast, the first integral, which involves $\phi_1$, grows without bound. With $w_0=c_1\phi_1$, the second term of $F_2$ contributes $c_1 d_\beta x^{-2/3}\phi_1^2$ to the integrand $\phi_1F_2$.
Since $\phi_1^2\sim(2/\pi)\cos^2\Phi$ oscillates about a non-zero mean,
$\langle\phi_1^2\rangle\to1/\pi$, this term has the non-oscillatory part
$c_1d_\beta/(\pi x^{2/3})$, which is not integrable at infinity. The first term of the source
has a non-zero mean as well, but carries the weight $x^{-4/3}$, which is integrable at infinity and contributes only a
constant. Hence

\begin{equation}
  \int_0^x \phi_1 F_2\,dy \;\sim\; \frac{3c_1d_\beta}{\pi}\,x^{1/3}\qquad(x\to\infty).
\label{eq:secular}
\end{equation}
Hence $w_2$ carries a slowly growing sine component,
\begin{equation}
  w_2 \sim c_1\sqrt{\tfrac{2}{\pi}}\Big[\epsilon_{2\beta}\cos\Phi
   + \tfrac{3d_\beta}{2}\,x^{1/3}\sin\Phi + \cdots\Big],
\label{eq:w2Large}
\end{equation}
which would apparently produce an $x$-dependent, potentially diverging contribution to the
amplitude. This is not the case, and is an artefact of truncating the inner solution at second
order: the growing sine shifts the phase of the oscillation rather than raising its amplitude.
Indeed, the inner equation \eqref{eq:wallEq} reads $w''+Qw=0$ with $Q$ as in \eqref{eq:normalform}, whose
WKB amplitude and phase are $Q^{-1/4}$ and $\int\!\sqrt{Q}\,dx$, respectively. The $\beta$-dependent corrections
to $Q$ decay in $x$, so $Q^{-1/4}$ does not grow, while the phase contains the term
$-\tfrac32 d_\beta\Lambda_n^{-2/3}x^{1/3}=O(\Lambda_n^{-1/2}\eta^{1/2})$ already present in
the matching phase \eqref{eq:largeinnerphase}. The growing sine in \eqref{eq:w2Large} is the
linearization of this phase term, and the norm, which depends on the bounded amplitude
$Q^{-1/4}$, is unaffected at the retained order. 

Collecting the three orders, the inner solution
$w=w_0+\Lambda_n^{-1/3}w_1+\Lambda_n^{-2/3}w_2$ has a large-$x$ cosine component with amplitude
$1+\epsilon_\beta\Lambda_n^{-1/3}+\epsilon_{2\beta}\Lambda_n^{-2/3}$ (from $w_0$, $w_1$ and
$w_2$ respectively) and a sine component with amplitude $\alpha_\beta\Lambda_n^{-1/3}$ (from
$w_1$), whose combined amplitude
$$
\sqrt{\left(1+\epsilon_\beta\Lambda_n^{-1/3}+\epsilon_{2\beta}\Lambda_n^{-2/3}\right)^2+\left(\alpha_\beta\Lambda_n^{-1/3}\right)^2}\simeq 1+\epsilon_\beta\Lambda_n^{-1/3}
+\Big(\epsilon_{2\beta}+\tfrac12\alpha_\beta^2\Big)\Lambda_n^{-2/3}+O(\Lambda_n^{-1})
$$
yields \eqref{2ndorderamplitude}.


\section{Eigenvalues and norms for $\beta=1/2$}
\label{appB}


\begin{table}[H]
  \centering
  \begin{minipage}{0.46\linewidth}
    \centering
    \begin{tabular}{c c c}
      \toprule
      $n$ & $\lambda_n$ & $C_n$ \\
      \midrule
       1 &  3.09145 & 2.065367 \\
       2 &  4.95947 & 1.010753 \\
       3 &  6.86327 & 0.690562 \\
       4 &  8.78382 & 0.533465 \\
       5 & 10.71476 & 0.439272 \\
       6 & 12.65295 & 0.376077 \\
       7 & 14.59655 & 0.330511 \\
       8 & 16.54439 & 0.295966 \\
       9 & 18.49565 & 0.268790 \\
      10 & 20.44976 & 0.246798 \\
      11 & 22.40627 & 0.228597 \\
      12 & 24.36486 & 0.213259 \\
      13 & 26.32525 & 0.200138 \\
      14 & 28.28724 & 0.188770 \\
      15 & 30.25064 & 0.178816 \\
      \bottomrule
    \end{tabular}
  \end{minipage}\hfill
  \begin{minipage}{0.46\linewidth}
    \centering
    \begin{tabular}{c c c}
      \toprule
      $n$ & $\lambda_n$ & $C_n$ \\
      \midrule
      16 & 32.21532 & 0.170018 \\
      17 & 34.18115 & 0.162180 \\
      18 & 36.14803 & 0.155146 \\
      19 & 38.11586 & 0.148795 \\
      20 & 40.08458 & 0.143028 \\
      21 & 42.05410 & 0.137765 \\
      22 & 44.02439 & 0.132940 \\
      23 & 45.99537 & 0.128499 \\
      24 & 47.96701 & 0.124396 \\
      25 & 49.93926 & 0.120592 \\
      26 & 51.91208 & 0.117054 \\
      27 & 53.88545 & 0.113754 \\
      28 & 55.85932 & 0.110668 \\
      29 & 57.83368 & 0.107775 \\
      30 & 59.80849 & 0.105057 \\
      \bottomrule
    \end{tabular}
  \end{minipage}
  \caption{First thirty eigenvalues $\lambda_n$ and norms $C_n$ of the
    Chen--Libby problem for $\beta=1/2$, computed by two-sided shooting. Norms are in the convention
    $C_n=\int_0^\infty (f')^4 z_n^2\,\mathrm{d}\eta$ with $N_n''(0)=1$.}
  \label{tab:thirty_modes}
\end{table}

\end{appendices}

\end{document}